\documentclass[10pt,a4paper]{article}

\usepackage[final]{optional}

\usepackage{amsfonts, amsmath, wasysym}

\usepackage{tikz}

\usepackage{amssymb,amsthm,
paralist
}

\usepackage{
latexsym,
}

\usepackage{url}

\definecolor{darkgreen}{rgb}{0,0.5,0}
\definecolor{darkred}{rgb}{0.7,0,0}
\usepackage[colorlinks, 
citecolor=darkgreen, linkcolor=darkred
]{hyperref}

\theoremstyle{plain}
\newtheorem{lemma}{Lemma}[section]
\newtheorem{thm}[lemma]{Theorem}

\theoremstyle{definition}

\newtheorem{problem}[lemma]{Problem}

\newtheorem{rmk}[lemma]{Remark}

\numberwithin{equation}{section}

\newcommand{\pl}[2]{{\frac{\partial #1}{\partial #2}}}

\newcommand{\be}{\beta}

\newcommand{\de}{\delta}

\newcommand{\la}{\lambda}
\newcommand{\La}{\Lambda}

\newcommand{\Si}{\Sigma}

\newcommand{\ep}{\varepsilon}

\newcommand{\R}{\ensuremath{{\mathbb R}}}
\newcommand{\N}{\ensuremath{{\mathbb N}}}

\newcommand{\Z}{\ensuremath{{\mathbb Z}}}

\newcommand{\downto}{\downarrow}

\newcommand{\lap}{\Delta}

\DeclareMathOperator{\Area}{Area}

\def\blbox{\quad \vrule height7.5pt width4.17pt depth0pt}

\newcommand{\beq}{\begin{equation}}
\newcommand{\eeq}{\end{equation}}
\newcommand{\beqa}{\begin{equation}\begin{aligned}}
\newcommand{\eeqa}{\end{aligned}\end{equation}}
\newcommand{\brmk}{\begin{rmk}}
\newcommand{\ermk}{\end{rmk}}
\newcommand{\partref}[1]{\hbox{(\csname @roman\endcsname{\ref{#1}})}}
\newcommand{\half}{\frac{1}{2}}

\newcommand{\cmt}[1]{\opt{draft}{\textcolor[rgb]{0.5,0,0}{
$\LHD$ #1 $\RHD$\marginpar{\blbox}}}}

\newcommand{\Rm}{{\mathrm{Rm}}}

\usepackage{soul}

\title{{
\bf
Loss of initial data under \\ limits of Ricci flows
\footnote{
This is the submitted version from 2019, but with updated references. Problem \ref{the_prob} can be answered using the theory in \cite{TY3}.
Appeared in 
`Minimal surfaces: integrable systems and visualisation,' 
[Granada, 2018] 
T. Hoffmann, M. Kilian, K. Leschke, F. Martin (Eds.) 
Springer Proceedings in Mathematics \& Statistics \textbf{349}, 2021.
}
} 
\\ 
\cmt{DRAFT with comments}
}
\author{Peter M. Topping}
\date{1 September 2021}

\begin{document}

\maketitle
\parskip=10pt

\begin{abstract}
We construct a sequence of smooth Ricci flows on $T^2$, with standard uniform $C/t$ curvature decay, and with initial metrics converging to the standard flat unit-area square torus $g_0$ in the Gromov-Hausdorff sense, with the property that the flows themselves converge not to the static Ricci flow $g(t)\equiv g_0$, but to the static Ricci flow $g(t)\equiv 2g_0$ of twice the area.
\end{abstract}


\section{Introduction}

When tasked with starting a Ricci flow with singular initial data, the standard approach is to approximate the singular initial data by smooth initial Riemannian metrics, then to flow each of the smooth metrics and take a smooth limit of the smooth resulting flows. 
As an example, in \cite{MT18, ST2, hochard, TR, Miles3D} one flows so-called Ricci limit spaces, obtained as non-collapsed Gromov-Hausdorff limits 
$(X,d_X)$ of sequences of smooth 3-manifolds $(M_i,g_i)$ satisfying uniform lower Ricci bounds.
Each $(M_i,g_i)$ gives rise to a Ricci flow $(M_i,g_i(t))$ 
of one form or another,
with uniform curvature bounds $|\Rm|_{g_i(t)}\leq C/t$ for $t\in (0,T)$, and  Hamilton-Cheeger-Gromov compactness allows one to extract a subsequence that converges to a smooth limit Ricci flow $(M,g(t))$ for $t\in (0,T)$.

The challenge then is to show that the desired initial data is not lost in the limit $i\to\infty$. In other words, we require that the smooth limit Ricci flow $(M,g(t))$ converges weakly to the desired initial data $(X,d_X)$ as $t\downto 0$. This amounts to showing that the Riemannian distance $d_{g(t)}$ has a uniform limit $d_0$ as $t\downto 0$, and that  $(M,d_0)$ is a metric space that is isometric to $(X,d_X)$.
In \cite{ST2, Miles3D}, this was achieved by proving uniform lower Ricci bounds on the flows $g_i(t)$. In particular, the so-called double bootstrap technique of $\cite{ST1}$ gives the local lower Ricci control that implies the necessary control on the evolution of distances, and higher-dimensional versions in the presence of stronger curvature hypotheses can be found in \cite{BCW, YL}, with the closest analogue (also being purely local) in \cite{Hoch_thesis}.

In this note we clarify that without uniform lower Ricci bounds this programme will fail completely in general. Indeed, this loss of initial data can occur even when $(X,d_X)$ is a smooth manifold of any dimension $n\geq 2$ and the convergence is much stronger than Gromov-Hausdorff.

\begin{thm}
\label{supersize}
Let $(T^2,g_0)$ be the standard flat square torus arising as the quotient $\R^2/ \Z^2$, and 
let $d_0:T^2\times T^2\to [0,\infty)$ be the Riemannian distance corresponding to $g_0$.
Then there exists a sequence of smooth Ricci flows $g_i(t)$, $t\in [0,\infty)$ on $T^2$, with the property that $|K_{g_i(t)}|\leq c_0/t$ for some uniform $c_0$, 
all $i\in\N$, and all $t\in (0,\infty)$, and so that
\begin{enumerate}
\item
$d_{g_i(0)}\to d_0$ uniformly as $i\to\infty$, but 
\item
the Ricci flows $g_i(t)$ converge smoothly locally on $T^2\times (0,\infty)$ to the flat metric $2g_0$ of twice the volume.
\end{enumerate}
\end{thm}
Here $K_g$ denotes the Gauss curvature of a metric $g$.

Variations on Theorem \ref{supersize} show that it is not even necessary for a limit Ricci flow and a limit initial metric to share the same conformal structure. What is most important is that the limit Ricci flow should just be \emph{larger} than the limit initial metric.


Estimates of Hamilton-Perelman \cite{formations, P1} (see \cite{ST1} for the local version) tell us that the uniform $c_0/t$ decay on the curvature implies  lower semicontinuity at $t=0$ in the sense that 
$$d_{g(t)}(x,y)\geq d_0(x,y)-\be\sqrt{c_0 t},$$
for some universal $\be<\infty$ and for $t\geq 0$. Clearly there is no
analogous upper semicontinuity in this example, contrary to the situation in which there is a uniform lower Ricci bound.

In situations where we do have uniform estimates relating $d_0$ or $d_{g_i(0)}$
with $d_{g_i(t)}$ for $t>0$,
one might expect better behaviour. Indeed, one might be able to analyse the limit alone if one had a positive answer to the following, cf. \cite{TR, alix}:


\begin{problem}
\label{the_prob}
Suppose $(M,g_0)$ is a smooth compact Riemannian manifold and $g(t)$ is a smooth 
Ricci flow on $M$ for $t\in (0,T)$ with the property that $d_{g(t)}\to d_{g_0}$ uniformly as $t\downto 0$. Is it then true that $g(t)$ extends to a smooth Ricci flow for
$t\in [0,T)$ with $g(0)=g_0$?
\end{problem}

Thus the question is whether attainment of initial data in a metric sense implies attainment of the initial data smoothly. 

This problem is open even in extremely simple situations such as when 
$(M,g_0)$ is the flat unit square torus as before. The problem is then to show that the only Ricci flow $g(t)$ on $M$ for $t\in (0,\ep)$ with $d_{g(t)}\to d_{g_0}$ uniformly as $t\downto 0$ is the stationary flow $g(t)\equiv g_0$.
It is not even immediate that such a Ricci flow is conformally equivalent to $g_0$.
What \emph{is} currently understood in this situation, thanks to the work of T. Richard \cite{TR}, is that if we additionally impose a hypothesis of a lower Ricci bound for $g(t)$ (equivalent here to a lower Gauss curvature bound) then we must indeed have $g(t)\equiv g_0$.
The higher dimensional case is addressed in concurrent work of A. Deruelle, F. Schulze and M. Simon \cite{DSS} in the presence of both an upper $c_0/t$ curvature bound and a uniform lower Ricci bound.

\section{The construction}

With respect to coordinates on $T^2$ coming from the Euclidean coordinates $x,y$, we can write $g_0=dx^2+dy^2$. 
Given an initial metric $u_0g_0$, where $u_0:T^2\to (0,\infty)$ is smooth, there exists a  unique subsequent Ricci flow $g(t)$  of the form $u\/g_0$, where the smooth function $u:T^2\times [0,\infty)\to (0,\infty)$ solves
\begin{equation*}
\left\{
\begin{aligned}
\pl{u}{t}=\lap\log u & \qquad  \text{ on }T^2\times (0,\infty)\\
u(\cdot,0)=u_0 & \qquad  \text{ on } T^2,
\end{aligned}
\right.
\end{equation*}
as described by Hamilton \cite{Ham88} (see \cite{GT2, ICRF_UNIQ} for the general theory on surfaces). By Gauss-Bonnet, the area is constant, because
$$\frac{d}{dt}\Area (T^2,g(t))=-2\int K\,dV=0,$$
see e.g. \cite[(2.5.8)]{RFnotes}.
It is often convenient to work with the function $v:=\half\log u$, which then satisfies the equation $\pl{v}{t}=e^{-2v}\lap v$, and we view this solution lifted to $\R^2$ whenever convenient. A solution $v$ whose norm is initially bounded by $\Lambda$ retains this bound by the maximum principle, and parabolic regularity theory gives $C^k$ bounds at, say, time $t=1$, depending only on $k$ and $\Lambda$. 
(See, for example, the discussion in Appendix B of \cite{GGthesis}.)
In particular, by applying this estimate to rescaled solutions
$(x,t)\mapsto v(\la x,\la^2 t)$, for $\la>0$, we find that 
$$|D^k v|(\cdot,t)\leq \frac{C(k,\La)}{t^{k/2}},$$ 
and in particular we can control the Gauss curvature $K=-e^{-2v}\lap v$ by
\beq
|K_{g(t)}|\leq \frac{c_0(\La)}{t},
\eeq
(In fact, we always have $K_{g(t)}\geq -\frac{1}{2t}$, see e.g. 
\cite[Corollary 3.2.5]{RFnotes}.)

The specific Ricci flows $g_i(t)$ will be determined by their initial data $g_i(0)$, which in turn will be chosen to be an appropriate $i$-dependent function times $g_0$. 
Thus as above we can  write $g_i(t)=u_i(t)g_0$,
for a one parameter family of functions $u_i(t):T^2\to (0,\infty)$. Our task is to choose the functions $u_i(0)$ appropriately.

We will choose the functions $u_i(0)$ to lie always within $[1,2]$. As above, this property is then preserved by the flow, i.e. $u_i(t)\in [1,2]$ throughout $T^2$ and for each $t\geq 0$.

For each $i$, consider the lattice $L$ of points in $T^2$ represented by points 
$(a/i,b/i)$ in $\R^2$, where $a,b\in \{0,1,\ldots, i-1\}$. For each pair of points in this lattice, choose a minimising geodesic connecting them within $(T^2,g_0)$. 
Denote the union of the images of this finite number of geodesics by $\Si\subset T^2$. 
We ask that $u_i(0)$ takes the value $1$ on the whole of $\Si$. We can then extend to a smooth function $u_i(0):T^2\to [1,2]$ with almost-maximal area in the sense that $\Area(T^2,g_i(0))\geq 2-1/i$. (Note that the area would be exactly $2$ if we could choose $u_i(0)\equiv 2$.) In particular, $\|u_i(0)-2\|_{L^1(T^2,g_0)}\leq 1/i$.


Since $g_i(0)\geq g_0$,  the distance $d_{g_i(0)}(p,q)$ between any two points $p,q\in T^2$ with respect to $g_i(0)$ must be at least  $d_{0}(p,q)$.
On the other hand, we can always find a point $P$ in the lattice $L$
such that
the distance from $p$ to $P$
is less than $1/i$ when measured with respect to $g_0$ or even with respect to $2g_0$ or $g_i(0)$. 
Similarly we can find a lattice point $Q$ close to $q$. Since $u_i(0)=1$ on $\Si$, the distance between lattice points with respect to $g_i(0)$ is equal to the  distance with respect to $g_0$. Thus 
\beqa
d_{g_i(0)}(p,q)&\leq d_{g_i(0)}(p,P)+d_{g_i(0)}(P,Q)+d_{g_i(0)}(Q,q)\\
&\leq d_{g_0}(P,Q)+\frac2i\\
&\leq d_{g_0}(P,p)+d_{g_0}(p,q)+d_{g_0}(q,Q)+\frac2i\\
&\leq d_{g_0}(p,q)+\frac4i
\eeqa
and we see that $d_{g_i(0)}(p,q)$ converges uniformly to $d_0(p,q)$ as $i\to\infty$, as required.

We now turn to the subsequent flows $g_i(t)=u_i(t)g_0$. By the discussion above,
we have $\Area (T^2,g_i(t))=\Area (T^2,g_i(0))\geq 2-1/i$, or equivalently
\beq
\label{lim_to_2}
\|u_i(t)-2\|_{L^1(T^2,g_0)}\leq 1/i.
\eeq
Moreover, the flows $g_i(t)$ satisfy a uniform Gauss curvature estimate $|K_{g_i(t)}|\leq \frac{c_0}{t}$, for some universal $c_0$, as required. Their conformal factors also enjoy uniform $C^k$ bounds for any $k\in\N$ over $T^2\times [\de,\infty)$, any $\de>0$, where the bounds depend on $k$ and $\de$. Thus a subsequence will converge smoothly locally  on $T^2\times (0,\infty)$
(as tensors) to a limit Ricci flow $g(t)=u(t)g_0$. 
By passing \eqref{lim_to_2} to the limit, we find that $g(t)\equiv 2g_0$.

{\sc mathematics institute, university of warwick, coventry, CV4 7AL,
uk}\\
\url{https://www.warwick.ac.uk/~maseq}
\end{document}